# Hyperelliptic curves, continued fractions, and Somos sequences

**Alfred J. van der Poorten**[1,*]

*Centre for Number Theory Research, Sydney*

**Abstract:** We detail the continued fraction expansion of the square root of a monic polynomials of even degree. We note that each step of the expansion corresponds to addition of the divisor at infinity, and interpret the data yielded by the general expansion. In the quartic and sextic cases we observe explicitly that the parameters appearing in the continued fraction expansion yield integer sequences defined by bilinear relations instancing sequences of Somos type.

The sequence $\ldots, 3, 2, 1, 1, 1, 1, 1, 2, 3, 5, 11, 37, 83, \ldots$ is produced by the recursive definition

$$B_{h+3} = (B_{h-1}B_{h+2} + B_hB_{h+1})/B_{h-2} \qquad (1)$$

and consists entirely of integers.

On this matter Don Zagier comments [26] that 'the proof comes from the theory of elliptic curves, and can be expressed either in terms of the denominators of the co-ordinates of the multiples of a particular point on a particular elliptic curve, or in terms of special values of certain Jacobi theta functions.'

Below, I detail the continued fraction expansion of the square root of a monic polynomials of even degree. In the quartic and sextic cases I observe explicitly that the parameters appearing in the expansion yield integer sequences defined by relations including and generalising that of the example (1).

However it is well known, see for example Adams and Razar [1], that each step of the continued fraction expansion corresponds to addition of the divisor at infinity on the relevant elliptic or hyperelliptic curve; that readily explains Zagier's explanation.

## 1. Some Brief Reminders

### *1.1. The numerical case*

We need little more than the following. Suppose $\omega$ is a quadratic irrational integer $\omega$, defined by $\omega^2 - t\omega + n = 0$, and greater than the other root $\overline{\omega}$ of its defining equation. Denote by $a$ the integer part of $\omega$ and set $P_0 = a - t$, $Q_0 = 1$. Then

[1]Centre for Number Theory Research, 1 Bimbil Place, Killara, Sydney, NSW 2071, Australia, e-mail: alf@math.mq.edu.au

[*]This survey was written at Brown University, Providence, Rhode Island where the author held the position of Mathematics Distinguished Visiting Professor, Spring semester, 2005. The author was also supported by his wife and by a grant from the Australian Research Council.

*AMS 2000 subject classifications:* primary 11A55, 11G05; secondary 14H05, 14H52.

*Keywords and phrases:* continued fraction expansion, function field of characteristic zero, hyperelliptic curve, Somos sequence.

the continued fraction expansion of $\omega_0 := (\omega + P_0)/Q_0$ is a two-sided sequence of lines, $h$ in $\mathbb{Z}$,

$$\frac{\omega + P_h}{Q_h} = a_h - \frac{\overline{\omega} + P_{h+1}}{Q_h}\,; \quad \text{in brief } \omega_h = a_h - \overline{\rho}_h\,,$$

with $(\omega + P_{h+1})(\overline{\omega} + P_{h+1}) = -Q_h Q_{h+1}$ defining the integer sequences $(P_h)$ and $(Q_h)$. Obviously $Q_0$ divides $(\omega + P_0)(\overline{\omega} + P_0)$. This suffices to ensure that the integrality of the sequence $(a_h)$ of *partial quotients* guarantees that always $Q_h$ divides the norm $n + tP_h + P_h^2 = (\omega + P_h)(\overline{\omega} + P_h)$.

**Comment 1.** Consider the $\mathbb{Z}$-module $\mathfrak{i}_h = \langle Q_h, P_h + \omega\rangle$. It is a pleasant exercise to confirm that $\mathfrak{i}_h$ is an ideal of the domain $\mathbb{Z}[\omega]$ if and only if $Q_h$ does indeed divide the norm $(\omega + P_h)(\overline{\omega} + P_h)$ of its numerator.

One says that a real quadratic irrational $\omega_h$ is *reduced* if

$$\omega_h > 1 \quad \text{and its conjugate } \overline{\omega}_h \text{ asatisfies } -1 < \overline{\omega}_h < 0\,.$$

If the *partial quotient* $a_h$ is always chosen as the integer part of $\omega_h$ then $\omega_0$ reduced entails all the $\omega_h$ and $\rho_h$ are reduced; and, remarkably, $a_h$ — which starts life as the integer part of $\omega_h$ — always also is the integer part of $\rho_h$. Then conjugation of the continued fraction tableau retrieves the negative half of the expansion of $\omega_0$ from the continued fraction expansion of $\rho_0$.

**Comment 2.** What a continued fraction expansion does. Suppose $\alpha = [\,a_0\,, a_1\,, a_2\,, \ldots\,]$ with $\alpha > 1$, is our very favourite expansion, so much so that sometimes we go quite alpha — expanding arbitrary complex numbers $\beta = \beta_0$ by the 'alpha' rule

$$\beta_h = a_h + (\beta_h - a_h) \quad \text{and } (\beta_h - a_h)^{-1} = \beta_{h+1} = a_{h+1} + \cdots \text{ etc.}$$

What can one say about those 'alpha d' complete quotients $\beta_h$ ?

Quite a while ago, in 1836, Vincent reports that either (i) all $\beta_h > 1$, in which case $\beta = \alpha$; or (ii) for *all* sufficiently large $h$, $|\beta_h| < 1$ and the real part $\Re\beta_h$ of $\beta_h$ satisfies $-1 < \Re\beta_h < 0$; in other words, all those $\beta_h$ lie in the left hand half of the unit circle. *Proof*: Straightforward exercise; or see [4] or, better, the survey [2].

Of this result Uspensky [24] writes: 'This remarkable theorem was published by Vincent in 1836 in the first issue of *Liouville's Journal*, but later [was] so completely forgotten that no mention of it is found even in such a capital work as the *Enzyclopädie der mathematischen Wissenschaften.* Yet Vincent's theorem is the basis of the very efficient method for separating real roots ... '.

The bottom line is this: When we expand a real quadratic irrational $\alpha$ then, willy-nilly, by conjugation we also expand $\overline{\alpha}$. By Vincent's theorem, *its* complete quotients eventually arrive in the left hand-half of the unit circle and, once 'reduced', they stay that way.

One readily confirms that the integers $P_h$ and $Q_h$ are bounded by

$$0 < 2P_h + t < \omega - \overline{\omega} \quad \text{and} \quad 0 < Q_h < \omega - \overline{\omega}\,.$$

It follows by the box principle that the continued fraction expansion of $\omega$ is periodic. More, the adjustment whereby we replace $\omega$ by $\omega_0 = \omega + a - t$ arranges that $\omega_0$ is reduced. Yet more, by conjugating the tableau one sees immediately that (an observation credited to Galois) for any $h$ the expansion of $\omega_h$ is purely periodic.

**Comment 3.** On conjugating the tableau, a putative preperiod becomes a 'post-period' in the expansion of $\rho_h$; which is absurd.

### 1.2. The function field case

Here, 'polynomial' replaces 'integer'. Specifically, set $Y^2 = D(X)$ where $D$ is a monic polynomial of even degree $\deg D = 2g+2$ and defined over the base field, $\mathbb{K}$ say. Then we may write

$$D(X) = A(X)^2 + 4R(X) \tag{2}$$

where $A$ is the polynomial part of the square root $Y$ of $D$, so $\deg A = g+1$, and the *remainder* $R$ satisfies $\deg R \le g$. It is in fact appropriate to study the expansion of $Z := \frac{1}{2}(Y + A)$. Plainly

$$\mathcal{C} : Z^2 - AZ - R = 0 \quad \text{with } \deg Z = g+1 \text{ and } \deg \overline{Z} < 0\,. \tag{3}$$

**Comment 4.** Note that $Y$ is given by a Laurent series $A + d_{-1}X^{-1} + d_{-2}X^{-2} + \cdots$ , an element of $\mathbb{K}((X^{-1}))$; in effect we expand around infinity. Felicitously, by restricting our attention to $Z$, and forgetting our opening remarks, the story we tell below makes sense over all base fields of arbitrary characteristic, including characteristic two. However, for convenience, below we mostly speak as if $\mathbb{K} = \mathbb{Q}$.

In the present context we study the continued fraction expansion of an element $Z_0$ of $\mathbb{K}(Z, X)$ leading to the expansion consisting of a tableau of lines, $h \in \mathbb{Z}$,

$$\frac{Z + P_h}{Q_h} = a_h - \frac{Z + P_{h+1}}{Q_{h+1}}\,, \quad \text{in brief } Z_h = a_h - \overline{R}_h\,, \text{ say,} \tag{4}$$

initiated by the conditions $\deg P_0 < g$, $\deg Q_0 \le g$ and $Q_0$ divides the norm $(Z + P_0)(\overline{Z} + P_0) = -R + P_0(A + P_0)$. Indeed, the story is *mutatis mautandis* precisely as in the numerical case, up to the fact that a function of $\mathbb{K}(Z, X)$ is reduced exactly when it has positive degree but its conjugate has negative degree. Here, analogously, we find that therefore all the $P_h$ and $Q_h$ satisfy

$$\deg P_h < g \quad \text{and} \quad \deg Q_h \le g\,, \tag{5}$$

the conditions equivalent to the $Z_h$ and $R_h$ all being reduced.

#### 1.2.1. Quasi-periodicity

If the base field $\mathbb{K}$ is infinite then the box principle does *not* entail periodicity. In a detailed reminder exposition on continued fractions in quadratic function fields at Section 4 of [17], we are reminded that periodicity entails the existence of a non-trivial unit, of degree $m$ say, in $\mathbb{K}[Z, X]$. Conversely however, the exceptional existence of such a unit implies only 'quasi-periodicity' — in effect, periodicity 'twisted' by multiplication of the period by a nonzero element of $\mathbb{K}$. The existence of an *exceptional* unit entails the divisor at infinity on the curve $\mathcal{C}$ being torsion of order dividing $m$. If quasi-periodic, the expansion of the reduced element $Z_0 = (Z + P_0)/Q_0$ is purely quasi-periodic.

**Comment 5.** Consider the surprising integral

$$\begin{aligned}\int \frac{6x\,dx}{\sqrt{x^4 + 4x^3 - 6x^2 + 4x + 1}} &= \log\Big(x^6 + 12x^5 + 45x^4 + 44x^3 - 33x^2 + 43 \\ &\qquad + (x^4 + 10x^3 + 30x^2 + 22x - 11)\sqrt{x^4 + 4x^3 - 6x^2 + 4x + 1}\Big),\end{aligned}$$

a nice example of a class of *pseudo-elliptic* integrals

$$\int \frac{f(x)dx}{\sqrt{D(x)}} = \log\bigl(a(x) + b(x)\sqrt{D(x)}\,\bigr). \tag{6}$$

Here we take $D$ to be a monic polynomial defined over $\mathbb{Q}$, of even degree $2g+2$, and not the square of a polynomial; $f$, $a$, and $b$ denote appropriate polynomials. We suppose $a$ to be nonzero, say of degree $m$ at least $g+1$. One sees readily that necessarily $\deg b = m - g - 1$, that $\deg f = g$, and that $f$ has leading coefficient $m$. In the example, $m = 6$ and $g = 1$.

The trick is to recognise that obviously $a^2 - b^2D$ is a non-zero constant and $a + b\sqrt{D}$ is a unit of degree $m$ in the domain $\mathbb{Q}(x, \sqrt{D(x)}$ and is not necessarily of norm $\pm 1$ — it is this that corresponds to *quasi*-periodicity; for details see [14].

### *1.2.2. Normaility of the expansion*

In the sequel, I suppose that $Z_0$ has been so chosen that its continued fraction expansion is *normal*: namely, all its partial quotients are of degree 1. This is the generic case if $\mathbb{K}$ is infinite. Since I have the case $\mathbb{K} = \mathbb{Q}$ in mind, I refer to elements of $\mathbb{K}$ as 'rational numbers'.

**Comment 6.** Naïvely, it is not quite obvious that the case all partial quotients of degree one is generic, let alone that this generic situation can be freely arranged in our partcular situation. I comment on the latter matter immediately below but the former point is this: when one inverts $d_{-1}X^{-1} + d_{-2}X^{-2} + \cdots$, obtaining a polynomial plus $e_{-1}X^{-1} + e_{-2}X^{-2} + \cdots$ it is highly improbable that $e_{-1}$ vanishes because $e_{-1}$ is a somewhat complicated rational function of several of the $d_{-i}$. See my remarks in [15].

### *1.2.3. Ideal classes*

If $\mathbb{K}$ is infinite, choosing $P_0$ and $Q_0$ is a matter of selecting one of infinitely many different ideal classes of $\mathbb{Z}$-modules $\{\langle Q_h, Z + P_h\rangle : h \in \mathbb{Z}\}$. Only a thin subset of such classes fails to give rise to a normal continued fraction expansion. Of course our choice of $P_0$ and $Q_0$ will certainly avoid the blatantly *singular* principal class: containing the ideal $\langle 1, Z\rangle$.

### *1.2.4. Addition on the Jacobian*

Here's what the continued fraction expansion does. The set of zeros $\{\omega_{h,1}, \ldots, \omega_{h,g}\}$ of $Q_h$ defines a rational divisor on the hyperelliptic curve $\mathcal{C}$. In plain language: any one of these zeros $\omega_h$ is the $X$ co-ordinate of a point on $\mathcal{C}$, here viewing $\mathcal{C}$ as defined over the algebraic extension $\mathbb{K}(\omega_h)$, generically of degree $g$ over the base field $\mathbb{K}$. In particular, in the case $g = 1$, when $\mathcal{C}$ is an elliptic curve, the unique zero $w_h$ of $Q_h$ provides a point on $\mathcal{C}$ defined over the base field; see §2.2 at page 217 below.

Selecting $P_0$ and $Q_0$ is to choose a divisor class, say $M$, thus a point on the Jacobian $\mathrm{Jac}(\mathcal{C})$ of the curve $\mathcal{C}$ defining $Z$. Let $S$ denote the divisor class defined by the divisor at infinity. Then each complete quotient $Z_h = (Z + P_h)/Q_h$ has divisor in the class $M_{h+1} := M + hS$. I show this in [17] for $g = 1$ making explicit remarks of Adams and Razar [1]. See further comment at §2.2. For higher genus cases, one wilkl find helpful the introduction to David Cantor's paper [7] and the instructive discussion by Kristin Lauter in [13]. A central theme of the paper [3] is a generalisation of the phenomenon to Padé approximation in arbitrary algebraic function fields. My suggestion that partial quotients of degree one are generic in our examples is the same remark as that divisors on $\mathcal{C}/\mathbb{K}$ typically are given by $g$ points on $\mathcal{C}/\mathbb{F}$, where $\mathbb{F}$ is some algebraic extension of $\mathbb{K}$.

## 2. The Continued Fraction Expansion

Evidently, the polynomials $P_h$ and $Q_h$ in

$$\frac{Z+P_h}{Q_h} = a_h - \frac{Z+P_{h+1}}{Q_{h+1}}\,, \tag{7}$$

are given sequentially by the formulas

$$\begin{aligned} &P_h + P_{h+1} + A = a_h Q_h \quad \text{and} \\ &-Q_h Q_{h+1} = (Z+P_{h+1})(\overline{Z}+P_{h+1}) = -R + P_{h+1}(A+P_{h+1}). \end{aligned} \tag{8}$$

### *2.1. A naïve approach*

At first glance one might well be tempted to use this data by spelling out the first recursion in terms of $g$ equations linearly relating the coeficients of the $P_h$, and the second recursion in terms of $2g+1$ equations quadratically relating the coefficients of the $Q_h$. Even for $g=1$, doing that leads to a fairly complicated analysis; see [16], and the rather less clumsy [17]. For $g=2$, I had to suffice myself with the special case $\deg R = 1$ so that several miracles (I called it 'a ridiculous computation') could yield a satisfying result [18]. I realised that a quite different view of the problem would be needed to say anything useful in higher genus cases.

### *2.2. Less explicit use of the recursions*

Denote by $\omega_h$ a zero of $Q_h(X)$. Then

$$\begin{aligned} &P_h(\omega_h) + P_{h+1}(\omega_h) + A(\omega_h) = 0 \\ &\qquad \text{so } -Q_h Q_{h+1} = -R + P_{h+1}\bigl(A+P_{h+1}\bigr) \\ &\qquad \text{becomes } R(\omega_h) = -P_{h+1}(\omega_h)P_h(\omega_h). \end{aligned} \tag{9}$$

This, together with a cute trick, already suffices to tame the $g=1$ case: that is, the case beginning as the square root of a monic quartic polynomial. Indeed, if $g=1$ then $\deg P_h = g-1$; so the $P_h$ are constants, say $P_h(X) =: e_h$. Also, $\deg Q_h = 1$, say $Q_h(X) =: v_h(X-w_h)$, so $\omega_h = w_h$. Further, $\deg R \le 1$, say $R(X) =: v(X-w)$.

First, (9) tells us directly that $v(w-w_h) = e_h e_{h+1}$. Second, this is the 'cute trick', obviously $-R(w_h) = v(w-w_h) = Q_h(w)\cdot v/v_h$. An intelligent glance at (8) reminds us that $v_{h-1}v_h = -e_h$. Hence, by (8) and $R(w)=0$,

$$e_{h-1}e_h^2 e_{h+1} = -Q_{h-1}(w)Q_h(w)\cdot v^2/e_h = v^2\bigl(e_h + A(w)\bigr)\,.$$

The exciting thing about the recursion

$$e_{h-1}e_h^2e_{h+1} = v^2\bigl(e_h + A(w)\bigr) \tag{10}$$

is that, among the parameters varying with $h$, it involves the $e_h$ alone; moreover, its coefficients $v^2$ and $v^2A(w)$ depend only on the curve $\mathcal{C}$ and not on the initial conditions $P_0 = e_0$ and $Q_0 = v_0(X - w_0)$.

More, if $X = w_h$ then $Z = -e_h$ or $Z = -e_{h+1}$; yielding pairs of rational points on $\mathcal{C}$. Moreover, the transformation $U = Z$, $V - v = XZ$ transforms the curve $\mathcal{C}$ to a familiar cubic model

$$\mathcal{E} : V^2 - vV = \text{monic cubic in } U \text{ with zero constant coefficient,} \tag{11}$$

essentially by moving one of the two points $S_{\mathcal{C}}$, say, at infinity on $\mathcal{C}$ to the origin $S_{\mathcal{E}} = (0,0)$ on $\mathcal{E}$. As at §1.2.4 on page 215 above, denote the point $(w_h, -e_h)$ on $\mathcal{C}$ by $M_{h+1}$. Recall that an elliptic curve is an abelian group with group operation denoted by $+$. Set $M_1 = M$. One confirms, see [17] for details, that $M_{h+1} = M + hS_{\mathcal{C}}$ by seeing that this plainly holds on $\mathcal{E}$ where the addition law is that $F + G + H = 0$ if the three points $F$, $G$, $H$ on $\mathcal{E}$ lie on a straight line.

I assert at §1.2.4 that precisely this property holds also in higher genus. There, however, one is forced to use the gobbledegook language of 'divisor classes on the Jacobian of the curve' in place of the innocent 'points on the curve' allowed in the elliptic case.

**Comment 7.** By the way, because each rational number $-e_h$ is the $U$ co-ordinate of a rational point on $\mathcal{E}$ it follows that its denominator must be the square of an integer.

### *2.3. More surely useful formulas*

Set $a_h(X) = (X + \nu_h)/u_h$; so $u_h$ is the leading coefficient of $Q_h$ (its coefficient of $X^g$) and $a_h$ vanishes at $-\nu_h$. Below, we presume that $d_h$ denotes the leading coefficient of $P_h$ (its coefficient of $X^{g-1}$). Then

$$\begin{aligned} R(X) + P_{h+1}(X)P_h(X) &= Q_h(X)Q_{h+1}(X) + P_{h+1}(X)(X+\nu_h)Q_h(X)/u_h \\ &= Q_{h-1}(X)Q_h(X) + P_h(X)(X+\nu_h)Q_h(X)/u_h\,. \end{aligned}$$

Plainly, we should divide by $Q_h(X)/u_h$ and may set

$$\begin{aligned} C_h(X) &:= \bigl(R(X) + P_{h+1}(X)P_h(X)\bigr)/\bigl(Q_h(X)/u_h\bigr) \\ &= u_hQ_{h+1}(X) + P_{h+1}(X)(X+\nu_h) \\ &= u_hQ_{h-1}(X) + P_h(X)(X+\nu_h). \end{aligned} \tag{12}$$

Since $\deg R \le g$, and the $P$ all have degree $g-1$, and the $Q$ degree $g$, it follows that the polynomial $C$ has degree $g-2$ if $g \ge 2$, or is constant in the case $g = 1$. If $g = 2$ it of course also is constant and then, with $R(X) = u(X^2 - vX + w)$, its leading coefficient is $d_hd_{h+1} + u$ so we have, identically, $C_h(X) = d_hd_{h+1} + u$. If $g \ge 3$ then $C_h$ is a polynomial with leading coefficient $d_hd_{h+1}$.

### *2.4. The case $g = 2$*

If $\deg A = 3$ then $R = u(X^2 - vX + w)$ is the most general remainder. In the continued fraction expansion we may take $P_h(X) = d_h(X + e_h)$. As just above, denote by $u_h$ the leading coefficient of $Q_h(X)$ and note that $u_{h-1}u_h = -d_h$. In the case $g = 1$ we succeeded in finding an identity in just the parameters $e_h$; here we seek an identity just in the $d_h$.

First we note, using (12) and $C_h(X) = d_h d_{h+1} + u$, that

$$d_h d_{h+1} + u = C_h(-e_h) = u_h R(-e_h)/Q_h(-e_h) = u_h Q_{h-1}(-e_h)\,. \tag{13}$$

$$d_{h-1}d_h + u = C_{h-1}(-e_h) = u_{h-1}R(-e_h)/Q_{h-1}(-e_h) = u_{h-1}Q_h(-e_h)\,. \tag{14}$$

Hence, cutely,

$$(d_{h-1}d_h + u)(d_h d_{h+1} + u) = u_{h-1}u_h R(-e_h) = -ud_h(e_h^2 + ve_h + w)\,. \tag{15}$$

Set $R(X) =: u(X-\omega)(X-\overline{\omega})$. Also from (12), we have

$$C_h(\omega)Q_h(\omega) = P_{h+1}(\omega)P_h(\omega)$$

and therefore

$$C_{h-1}(\omega)C_h(\omega)Q_{h-1}(\omega)Q_h(\omega) = P_{h-1}(\omega)\big(P_h(\omega)\big)^2 P_{h+1}(\omega)\,.$$

But $-Q_{h-1}(\omega)Q_h(\omega) = P_h(\omega)\big(A(\omega) + P_h(\omega)\big)$. Together with (15) we find that

$$\begin{aligned} &u^3\big(A(\omega) + d_h(\omega + e_h)\big)\big(A(\overline{\omega}) + d_h(\overline{\omega} + e_h)\big) \\ &\qquad = d_{h-1}d_h^3 d_{h+1}(d_{h-2}d_{h-1} + u)(d_{h+1}d_{h+2} + u). \end{aligned} \tag{16}$$

In the special case $u = 0$, that is: $R(X) = -v(X - w)$, a little less argument yields the more amenable

$$-v^3\big(A(\omega) + d_h(w + e_h)\big) = d_{h-2}d_{h-1}^2 d_h^3 d_{h+1}^2 d_{h+2}\,.$$

In this case, we have $-vd_h(w + e_h) = -d_h R(-e_h) = d_{h-1}d_h^2 d_{h+1}$, finally providing

$$d_{h-2}d_{h-1}^2 d_h^3 d_{h+1}^2 d_{h+2} = v^2 d_{h-1}d_h^2 d_{h+1} - v^3 A(w)\,. \tag{17}$$

**Comment 8.** This reasonably straightforward argument removes the 'miraculous' aspects from the corresponding discussion in [18]. Moreover, it gives *a* result for $u \neq 0$. However, I do not yet see how to remove the dependence on $e_h$ in (16) so as to obtain a polynomial relation in the $d$s.

## 3. Somos Sequences

The complexity of the various parameters in the continued fraction expansions increases at frantic pace with $h$. For instance the logarithm of the denominators of $e_h$ of the elliptic case at §2.2 is readily proved to be $O(h^2)$ and the same must therefore hold for the logarithmic height of each of the parameters.

Denote by $A_h^2$ the denominator of $e_h$ (recall the remark at Comment 7 on page 217 that these denominators in fact are squares of integers). It is a remarkable fact holding for the co-ordinates of multiples of a point on an elliptic curve that *in general*

$$A_{h-1}A_{h+1} = e_h A_h^2 \tag{18}$$

Moreover, it is a simple exercise to see that (18) entails $A_{h-2}A_{h+2} = e_{h-1}e_h^2 e_{h+1}A_h^2$. Thus, on multiplying the identity (10), namely $e_{h-1}e_h^2 e_{h+1} = v^2\big(e_h + A(w)\big)$, by $A_h^2$ we find that

$$A_{h-2}A_{h+2} = v^2 A_{h-1}A_{h+1} + v^2 A(w)A_h^2 \tag{19}$$

gives a quadratic recursion for the 'denominators' $A_h$.

**Comment 9.** My remark concerning the co-ordinates of multiples of a point on an elliptic curve is made explicit in Rachel Shipsey's thesis [21]. The fact this also holds after a translation is shown by Christine Swart [23] and is in effect again proved here by way of the recursion (10). My weaselling 'in general' is to avoid my having to chat on about exceptional primes — made evident by the equation defining the elliptic curve — at which the $A_h$ may not be integral at all. In other words, in true generality, there may be a finite set $\mathcal{S}$ of primes so that the $A_h$ actually are just $\mathcal{S}$-integers: that is, they may have denominators but primes dividing such denominators must belong to $\mathcal{S}$.

### *3.1. Michael Somos's sequences*

Some fifteen years ago, Michael Somos noticed [11, 20], that the two-sided sequence $C_{h-2}C_{h+2} = C_{h-1}C_{h+1} + C_h^2$, which I refer to as 4-Somos in his honour, apparently takes only integer values if we start from $C_{h-1} = C_h = C_{h+1} = C_{h+2} = 1$. Indeed Somos went on to investigate also the width 5 sequence, $B_{h-2}B_{h+3} = B_{h-1}B_{h+2} + B_hB_{h+1}$, now with five initial 1s, the width 6 sequence $D_{h-3}D_{h+3} = D_{h-2}D_{h+2} + D_{h-1}D_{h+1} + D_h^2$, and so on, testing whether each when initiated by an appropriate number of 1s yields only integers. Naturally, he asks: "What is going on here?"

By the way, while 4-Somos (A006720) 5-Somos (A006721), 6-Somos (A006722), and 7-Somos (A006723) all do yield only integers; 8-Somos does not. The codes in parentheses refer to Neil Sloane's 'On-line encyclopedia of integer sequences' [22].

### *3.2. Elliptic divisibility sequences*

Sequences generalising those considered by Somos were known in the literature. Morgan Ward had studied anti-symmetric sequences $(W_h)$ satisfying relations

$$W_{h-m}W_{h+m}W_n^2 = W_{h-n}W_{h+n}W_m^2 - W_{m-n}W_{m+n}W_h^2\,. \tag{20}$$

He shows that if $W_1 = 1$ and $W_2 \big| W_4$ then $a \big| b$ implies that $W_a \big| W_b$; that is, the sequences become *divisibility sequences* (compare the Fibonacci numbers). For a brief introduction see Chapter 12 of [9].

There is a drama here. The recurrence relation

$$W_{h-2}W_{h+2} = W_2^2 W_{h-1}W_{h+1} - W_1W_3W_h^2\,,$$

and four nonzero initial values, already suffices to produce $(W_h)$. Thus (20) for all $m$ and $n$ is apparently entailed by its special case $n = 1$ and $m = 2$. The issue is whether the definition (20) is coherent.

One has to go deep into Ward's memoir [25] to find an uncompelling proof that there is in fact a solution sequence, namely one defined in terms of quotients of Weierstrass sigma functions. More to the point, given integers $W_1 = 1$, $W_2$, $W_3$, and $W_4$, there always is an associated elliptic curve. In our terms, there is a curve $\mathcal{C} : Z^2 - AZ - R = 0$ with $\deg A = 2$, $\deg R = 1$ and the sequence $(W_h)$ arises from the continued fraction expansion of $Z_1 = Z/(-R)$. I call $(W_h)$ the *singular* sequence because in that case $e_1 = 0$ — so the partial quotient $a_0(X)$ is not linear.

My 'translated' sequences $(A_h)$ were extensively studied by Christine Swart in her thesis [23].

### *3.3. Somos sequences*

It is natural to generalise Michael Somos's questions and to study recurrences of the kind he considers but with coefficients and initial values not all necessarily 1. Then our recurrence (19) is a general instance of a Somos 4 sequence, an easy computation confirms that $-v = W_2$ and $v^2A(w) = -W_3$, and given the recursion and four consecutive $A_h$ one can readily identify the curve $\mathcal{C}$ and the initial 'translation' $M = (w_0, -e_0)$.

**Comment 10.** One might worry (I did worry) that $A_{h-2}A_{h+2} = aA_{h-1}A_{h+1} + bA_h^2$ does not give a rationally defined elliptic curve if $a$ is not a square. No worries. One perfectly happily gets a quadratic twist by $a$ of a rationally defined curve.

For example, 4-Somos, the sequence $(C_h) = (\dots, 2, 1, 1, 1, 1, 2, 3, 7, \dots)$ with $C_{h-2}C_{h+2} = C_{h-1}C_{h+1} + C_h^2$ arises from

$$\mathcal{C} : Z^2 - (X^2 - 3)Z - (X - 2) = 0 \quad \text{with } M = (1, -1);$$

equivalently from $\mathcal{E} : V^2 - V = U^3 + 3U^2 + 2U$ with $M_{\mathcal{E}} = (-1, 1)$.

Christine Swart and I found a nice inductive proof [19] that if $(A_h)$ satisfies (19) then for all integers $m$ and $n$,

$$A_{h-m}A_{h+m}W_n^2 = W_m^2 A_{h-n}A_{h+n} - W_{m-n}W_{m+n}A_h^2 .$$

Our argument obviates any need for talk of transcendental functions and is purely algebraic.

It also is plain that $A_{h-1}A_{h+1} = e_hA_h^2$ yields $A_{h-1}A_{h+2} = e_he_{h+1}A_hA_{h+1}$ and $A_{h-2}A_{h+3} = e_{h-1}e_h^2e_{h+1}^2e_{h+2}A_hA_{h+1}$. However, although (10) directly entails

$$e_{h-1}e_h^2e_{h+1}^2e_{h+2} = -v^2A(w)e_he_{h+1} + \bigl(v^4 + 2wv^3A(w)\bigr),$$

it requires some effort to see this. Whatever, (19) eventually also gives

$$W_1W_2A_{h-m}A_{h+m+1} = W_mW_{m+1}A_{h-1}A_{h+2} - W_{m-1}W_{m+2}A_hA_{h+1} . \qquad (21)$$

For details see [17] and [19].

The case $m = 2$ of (21) includes all Somos 5 sequences. The sequence 5-Somos, $(B_h) = (\dots, 2, 1, 1, 1, 1, 1, 2, 3, 5, 11, \dots)$ with $B_{h-2}B_{h+3} = B_{h-1}B_{h+2} + B_hB_{h+1}$, arises from

$$Z^2 - (X^2 - 29)Z + 48(X + 5) = 0 \quad \text{with } M = (-3, -8);$$

equivalently from $\mathcal{E} : V^2 + UV + 6V = U^3 + 7U^2 + 12U$ with $M_{\mathcal{E}} = (-2, -2)$.

Actually, a Somos 5 sequence $(A_h)$ may also be viewed as a pair $(A_{2h})$ and $(A_{2h+1})$ of Somos 4 sequences coming from the same elliptic curve but with different translations (in fact differing by half its point $S$); see the discussion in [17].

### *3.4. Higher genus Somos sequences*

My purpose in studying the elliptic case was to be able to make impact on higher genus cases. That's been only partly achieved, what with little more than (17) at page 218 to show for the effort. I tame (17) by defining a sequence $(T_h)$, one hopes of integers, by way of $T_{h-1}T_{h+1} = d_hT_h^2$. We already know that then $T_{h-2}T_{h+2} =$ $d_{h-1}d_h^2d_{h+1}T_h^2$; also $T_{h-3}T_{h+3} = d_{h-2}d_{h-1}^2d_h^3d_{h+1}^2d_{h+2}T_h^2$ follows with only small extra effort. Thus (17) yields

$$T_{h-3}T_{h+3} = v^3T_{h-2}T_{h+2} - v^2A(w)T_h^2\,. \tag{22}$$

Then the sequence $(T_h) = (\ldots,\ 2,\ 1,\ 1,\ 1,\ 1,\ 1,\ 1,\ 2,\ 3,\ 4,\ 8,\ 17,\ 50,\ \ldots)$ satisfying $T_{h-3}T_{h+3} = T_{h-2}T_{h+2} + T_h^2$ is readily seen to derive from the genus 2 curve

$$\mathcal{C} : Z^2 - Z(X^3 - 4X + 1) - (X - 2) = 0\,. \tag{23}$$

A relevant piece of the associated continued fraction expansion is

$$\begin{aligned}
\frac{Z + 2X - 1}{X^2 - 1} &= X - \frac{\overline{Z} + X}{X^2 - 1}\\
\frac{Z + X}{-(X^2 - 2)} &= -X - \frac{\overline{Z} + X - 1}{-(X^2 - 2)}\\
Z_0 := \frac{Z + X - 1}{X^2 - X - 1} &= X + 1 - \frac{\overline{Z} + X - 1}{X^2 - X - 1}\\
\frac{Z + X - 1}{-(X^2 - 2)} &= -X - \frac{\overline{Z} + X}{-(X^2 - 2)}\\
\frac{Z + X}{X^2 - 1} &= X - \frac{\overline{Z} + 2X - 1}{X^2 - 1}\\
&\cdots
\end{aligned}$$

illustrating that $M$ is the divisor class defined by the pair of points $(\varphi, \overline{\varphi})$ and $(\overline{\varphi}, \varphi)$ — here, $\varphi$ is the golden ratio, a happenstance that I expect will please adherents to the cult of Fibonacci.

**Comment 11.** There does remain an issue here. Although (10) is concocted on the presumption that $e_h$ is never zero, it continues to make sense if some $e_h$ should vanish — for instance in the singular case. However, (17) is flat out false if $d_{h-1}d_hd_{h+1} = 0$. Indeed, a calm study of the argument yielding (17) sees us dividing by $v(e_h + w)$ at a critical point. These considerations together with Cantor's results [8] suggest that (22) should always be reported as multiplied by $T_h$ and, more to the point, that the general genus 2 relation when $u \neq 0$ will be cubic rather than quadratic.

## 4. Other Viewpoints

My emphasis here has been on continued fraction expansions producing sequences $M + hS$ of divisors — in effect the polynomials $Q_h$ — obtained by repeatedly adding a divisor $S$ to a starting 'translation' $M$. That viewpoint hints at arithmetic reasons for the integrality of the Somos sequences but does not do that altogether convincingly in genus greater than one.

### *4.1. The Laurent phenomenon*

As it happens, the integrality of the Somos sequences is largely a combinatorial phenomenon. In brief, as an application of their theory of *cluster algebras*, Fomin and Zelevinsky [10] prove results amply including the following. Suppose the sequence $(y_h)$ is defined by a recursion

$$y_{h+n}y_h = \alpha y_{h+r}y_{h+n-r} + \beta y_{h+s}y_{h+n-s} + \gamma y_{h+t}y_{h+n-t}\,,$$

with $0 < t < s < r \le \frac{1}{2}n$. Then the $y_h$ are Laurent polynomials* in the variables $y_0$, $y_1$, …, $y_{n-1}$ and with coefficients in the ring $\mathbb{Z}[\alpha,\beta,\gamma]$. That deals with all four term and three term quadratic recursions and thus with the cases Somos 4 to Somos 7. Rather more is true than may be suggested by the given example.

### 4.2. Dynamic methods

Suppose we start from a Somos 4 relation $A_{h-2}A_{h+2} = \alpha A_{h-1}A_{h+1} + \beta A_{h^2}$ and appropriate initial values $A_0$, $A_1$, $A_2$, $A_3$. Then one obtains rationals $e_h = A_{h-1}A_{h+1}/A_h^2$ satisfying the difference equation

$$e_{h+1} = \frac{1}{e_{h-1}e_h}\Bigl(\alpha + \frac{\beta}{e_h}\Bigr).$$

The point is that this equation has a first integral given by

$$J := J(e_{h-1},e_h) = e_{h-1}e_h + \alpha\Bigl(\frac{1}{e_{h-1}} + \frac{1}{e_h}\Bigr) + \frac{\beta}{e_{h-1}e_h} = J(e_h,e_{h+1})$$

and one can now construct an underlying Weierstrass elliptic function $\wp$. Indeed, the readily checked assertion that given $y \in \mathbb{C}$ there are constants $\alpha$ and $\beta$ so that

$$\bigl(\wp(x+y)-\wp(y)\bigr)\bigl(\wp(x)-\wp(y)\bigr)^2\bigl(\wp(x-y)-\wp(y)\bigr) = -\alpha\bigl(\wp(x)-\wp(y)\bigr)+\beta$$

reveals all; particularly that $\alpha = \wp'(y)^2$, $\beta = \wp'(y)^2(\wp(2y)-\wp(y))$.

Specifically, after fixing $x$ by $J \equiv \wp''(x)$ one sees that

$$-e_h = \wp(x+ny) - \wp(y)\,.$$

A program of this genre is elegantly carried out by Andy Hone in [12] for Somos 4 sequences. In [5], the ideas of [12] are shown to allow a Somos 8 recursion to be associated with adding a divisor on a genus 2 curve. Incidentally, that result coheres with the guess mooted at Comment 11 above that there always then is a cubic relation of width 6.

## References

[1] Adams, William W. and Razar, Michael J. (1980). Multiples of points on elliptic curves and continued fractions. *Proc. London Math. Soc.* **41**, 481–498. MR591651.

[2] Alesina, Alberto and Galuzzi, Massimo (1998). A new proof of Vincent's theorem. *Enseign. Math.* (2) **44**.3–4, 219–256. MR1659208; Addendum (1999). *Ibid.* **45**.3–4, 379–380. MR1742339.

[3] Bombieri, Enrico and Cohen, Paula B. (1997). Siegel's lemma, Padé approximations and Jacobians' (with an appendix by Umberto Zannier, and dedicated to Enzio De Giorgi). *Ann. Scuola Norm. Sup. Pisa* Cl. Sci. (4) **25**, 155–178. MR1655513.

*A Laurent polynomial in the variable $x$ is a polynomial in $x$ and $x^{-1}$.

[4] Bombieri, Enrico and van der Poorten, Alfred J. (1995). Continued fractions of algebraic numbers. In *Computational Algebra and Number Theory* (Sydney, 1992) Math. Appl., 325, Kluwer Acad. Publ., Dordrecht, 137–152. MR1344927.

[5] Braden, Harry W., Enolskii, Victor Z., and Hone, Andrew N. W. (2005). Bilinear recurrences and addition formulæ for hyperelliptic sigma functions'. 15pp: at http://www.arxiv.org/math.NT/0501162.

[6] Cassels, J. W. S. and Flynn, E. V. (1966). *Prolegomena to a Middlebrow Arithmetic of Curves of Genus* 2, London Mathematical Society Lecture Note Series, 230. Cambridge University Press, Cambridge, 1996. xiv+219 pp. MR1406090.

[7] Cantor, David G. (1987). Computing in the Jacobian of a hyperelliptic curve. *Math. Comp.* **48**.177, 95–101. MR866101.

[8] Cantor, David G. (1994). On the analogue of the division polynomials for hyperelliptic curves. *J. für Math.* (Crelle), **447**, 91–145. MR1263171.

[9] Everest, Graham, van der Poorten, Alf, Shparlinski, Igor, and Ward, Thomas (2003). *Recurrence Sequences*. Mathematical Surveys and Monographs 104, American Mathematical Society, xiv+318pp. MR1990179.

[10] Fomin, Sergey and Zelevinsky, Andrei (2002). The Laurent phenomenon. *Adv. in Appl. Math.*, **28**, 119–144. MR1888840.
Also 21pp: at http://www.arxiv.org/math.CO/0104241.

[11] Gale, David (1991). The strange and surprising saga of the Somos sequences. *The Mathematical Intelligencer* **13**.1 (1991), 40–42; Somos sequence update. *Ibid.* **13**.4, 49–50.

[12] Hone, A. N. W. (2005). Elliptic curves and quadratic recurrence sequences. *Bull. London Math. Soc.* **37**, 161–171. MRMR2119015

[13] Lauter, Kristin E. (2003). The equivalence of the geometric and algebraic group laws for Jacobians of genus 2 curves. *Topics in Algebraic and Noncommutative Geometry* (Luminy/Annapolis, MD, 2001), 165–171, *Contemp. Math.*, **324**, Amer. Math. Soc., Providence, RI. MR1986121.

[14] Pappalardi, Francesco and van der Poorten, Alfred J. (2004). Pseudo-elliptic integrals, units, and torsion.
12pp: at http://www.arxiv.org/math.NT/0403228.

[15] van der Poorten, Alfred J. (1998). Formal power series and their continued fraction expansion. In *Algorithmic Number Theory* (Proc. Third International Symposium, ANTS-III, Portland, Oregon, June 1998), Springer Lecture Notes in Computer Science **1423**, 358–371. MR1726084.

[16] van der Poorten, Alfred J. (2004). Periodic continued fractions and elliptic curves. In *High Primes and Misdemeanours*. Lectures in Honour of the 60th Birthday of Hugh Cowie Williams, Fields Institute Communications **42**, American Mathematical Society, 353–365. MR2076259.

[17] van der Poorten, Alfred J. (2005). Elliptic curves and continued fractions. 12pp: at http://arxiv.org/math.NT/0403225.

[18] van der Poorten, Alfred J. (2005). Curves of genus 2, continued fractions, and Somos sequences. 6pp: at http://arxiv.org/math.NT/0412372.

[19] van der Poorten, Alfred J. and Swart, Christine S. (2005). Recurrence relations for elliptic sequences: every Somos 4 is a Somos $k$.
7pp: at http://arxiv.org/math.NT/0412293.

[20] Propp, Jim. *The Somos Sequence Site.*
http://www.math.wisc.edu/~propp/somos.html.

[21] Shipsey, Rachel (2000). Elliptic divisibility sequences, PhD Thesis, Goldsmiths College, University of London.
http://homepages.gold.ac.uk/rachel/.

[22] Sloane, Neil *On-Line Encyclopedia of Integer Sequences.*
http://www.research.att.com/~njas/sequences/.

[23] Swart, Christine (2003). Elliptic curves and related sequences. PhD Thesis, Royal Holloway, University of London.

[24] Uspensky, J. V. (1948). *Theory of Equations*, McGraw-Hill Book Company.

[25] Ward, Morgan (1948). Memoir on elliptic divisibility sequences *Amer. J. Math.* **70**, 31–74. MR0023275.

[26] Zagier, Don (1966) Problems posed at the St Andrews Colloquium, Solutions, 5th day; see
http://www-groups.dcs.st-and.ac.uk/~john/Zagier/Problems.html.